\documentclass[conference]{IEEEtran}
\IEEEoverridecommandlockouts
\usepackage{cite}
\usepackage{amsmath,amssymb,amsfonts}
\usepackage{algorithmic}
\usepackage{graphicx}
\usepackage{textcomp}
\newcommand*{\QEDB}{\hfill\ensuremath{\square}}%
\newtheorem{Example}{Example}[section]
\newtheorem{problem}{Problem}[section]

\newtheorem{t1}{Theorem}[section]

\newtheorem{remark}{Remark}[section]
\newtheorem{lemma}{Lemma}[section]
\newtheorem{algorithm}{\underline{Algorithm}}[section]
\def\BibTeX{{\rm B\kern-.05em{\sc i\kern-.025em b}\kern-.08em
    T\kern-.1667em\lower.7ex\hbox{E}\kern-.125emX}}

\newcommand{\dd}{{\tt d}}
\begin{document}

\title{Structured real radius of controllability for higher order LTI systems\\
\thanks{The authors acknowledge the financial support from Science and Engineering
Research Board (SERB) India under the grant EMR/2015/002156. The first author of the paper acknowledges the financial support from
CSIR, HRDG, Govt. of India as a PhD student
with sanction letter number 09/081(1280)/2016-EMR-I.
A shorter version of this paper is accepted for publication in the proceedings of 27th Mediterranean Conference
on Control and Automation (MED '19), July 2019, Akko, ISRAEL}
}

\author{\IEEEauthorblockN{ Tanay Saha}
\IEEEauthorblockA{
\textit{Indian Institute of Technology}\\
Kharagpur, India \\
sahatanay50@gmail.com}
\and
\IEEEauthorblockN{ Swanand Khare}
\IEEEauthorblockA{
\textit{Indian Institute of Technology}\\
Kharagpur, India \\
srkhare@maths.iitkgp.ac.in}
\and
\IEEEauthorblockN{ Subashish Datta}
\IEEEauthorblockA{
\textit{Indian Institute of Technology}\\
Delhi, India \\
subashish@ee.iitd.ac.in }
}
\maketitle
\begin{abstract}
In this paper, we consider the problem of computing the nearest
uncontrollable (C-uncontrollable) system to a given higher order system. The distance to the nearest
uncontrollable system, also termed as the radius of controllability, is a good
measure of gauging the numerical robustness of the given system with respect
to controllability. Here, we invoke the equivalence of C-controllability of a higher order system with
full rank property of a certain Toeplitz structured matrix. This enables us to pose the problem of
computing the radius of controllability as equivalent to the problem of computing the nearest structured
low rank approximation of this Toeplitz structured matrix. Through several numerical examples and comparison
with the benchmark numerical problem, we illustrate that our approach works well.
\end{abstract}
\section{Introduction}
Many real life phenomena can be modeled as linear time-invariant (LTI) higher order systems as defined in the following equation:
\begin{equation}
\label{higher}
P_\dd \frac{d^\dd x(t)}{dt^\dd}+ P_{\dd-1} \frac{d^{\dd-1}x(t)}{dt^{\dd-1}}+\dots+ P_1 \frac{dx(t)}{dt}+ P_0 x(t)=b u(t)
\end{equation}
Here $P_0,P_1,\dots,P_{\dd-1},P_\dd \in \mathbb R^{N \times N},~b\in \mathbb R^{N \times M}$.
It should be noted that the case when $\dd=1$, the system in \eqref{higher} becomes an LTI descriptor system. In addition if
we consider the leading coefficient matrix $P_\dd$ to be non-singular, the system in \eqref{higher} can easily be
transformed into the usual state space system. However, in many real life applications, the matrix $P_\dd$ turns out to be a singular matrix
for $\dd \geq 1$.
(see for instance \cite{duan2010analysis,datta2014computation}).

Controllability of the system in \eqref{higher} plays a key role in the analysis of such systems.
However, the notion of controllability gives only a yes/no type of answer. It may so happen that small
perturbations to system parameters may render the system uncontrollable. In such a case, it is not sufficient to know
only whether the given system is controllable, but it is advisable to know how far this system is from being uncontrollable.
In other words, one may be interested in knowing the nearest uncontrollable system to the given system.
The distance to the nearest uncontrollable system, also called the radius of controllability, is a robust measure of controllability.
Several studies in the past have been conducted in the area of computation of the radius of controllability for
state space systems \cite{wicks1991computing,hu2004real,gu2006fast,khare}
as well as for descriptor systems and higher order systems \cite{son,son2012structured,simon}.

In \cite{esing1}, the author first introduced the notion of real radius of controllability for state space systems of the form
$\overset{.}{x}(t) = Ax(t) + Bu(t)$ and many researchers have followed this work as evident in the literature
(see for instance \cite{boley1986measuring,elsner1991algorithm,wicks1991computing,hu2004real,gu2006fast,khare} and references therein).
In \cite{elsner1991algorithm}, the radius of controllability for state space system was computed by minimizing the
minimum singular value function of a matrix pencil, which is not numerically efficient. A similar concept was implemented in
\cite{hu2004real,wicks1991computing} to compute the radius of controllability. In \cite{khare}, the problem of computing the
radius of controllability for state space system was shown to be equivalent to a problem of computing
the nearest Structured Low Rank Approximation (SLRA) of a Toeplitz structured matrix.

In \cite{simon}, a higher order linear time invariant system was considered. The radius of
controllability was computed using generalized singular value decomposition of a matrix pair, though finding generalized
singular value decomposition may not be numerically efficient in general. Basically, a generalized
real perturbation value function is minimized in the complex plane, which is an one dimensional optimization problem. The structured radius of
controllability of higher order system was also discussed in \cite{son}.

In this paper we will illustrate an approach to compute the radius of controllability which involves computation of the nearest SLRA of
a structured matrix. Due to some physical constraints, it may happen that some of the entries of $P_i (i=0,1,\dots,\dd)$ and $b$ are not allowed to
be perturbed while computing the radius of controllability. In such cases, the proposed approach is shown to be easily applied
to compute the nearest uncontrollable system. The proposed technique is computationally efficient as illustrated subsequently and avoids the
computation of
generalized singular value decomposition as done in \cite{simon}.
\subsection*{Criterion for Controllability of Higher order Systems}
Here, we state various types of controllability concepts introduced in the literature for LTI higher order, descriptor systems.
For a comprehensive study
of the descriptor systems the reader is referred to \cite{dai,duan2010analysis,lecture }.
Controllability for linear descriptor systems is analyzed by some graph-theoretic approach in \cite{Reinschke}.
The controllability of the system in \eqref{higher} is discussed in the following result from \cite{son}.
\begin{t1}
\label{th2}
The Higher order system in \eqref{higher} is $C^\dd$-controllable if and only if
\begin{enumerate}
\item rank$([P(s)~b])=N,$ for all finite $s\in \mathbb C$, and
\item rank$([P_\dd~b])=N$
\end{enumerate}
where, $P(s):=P_\dd s^\dd +\dots+P_1 s+P_0$.
\end{t1}
\subsection*{Descriptor systems}
Descriptor system is a particular case of higher order system with $\dd=1$. Consider the following descriptor system
\begin{equation}
\label{des}
E\overset{.}{z}(t)=Az(t)+Bu(t);~ E,A \in \mathbb R^{n \times n},~B\in \mathbb R^{n \times m}
\end{equation}
Controllability of the descriptor system in \eqref{des} plays an important
role in our paper. Now we state the theorem regarding the criterion of controllability of descriptor systems, which directly follows from
Theorem \ref{th2}.
\begin{t1}
\label{th3}
The descriptor system in \eqref{des} is C-controllable if and only if
\begin{enumerate}
\item rank$[sE-A ~~B]=n$, for all finite $s \in \mathbb C$, and
\item rank$[E~~B]=n$.
\end{enumerate}
\end{t1}
\section{Radius of controllability}
As discussed in the previous section, we now formally define the real radius of controllability, denoted as $R_c$, for the system in \eqref{higher}.
In order to simplify the notation we write $(P_\dd,\dots,P_1,P_0,b)$ is either C-controllable/uncontrollable to mean the system in \eqref{higher}
is C-controllable/uncontrollable respectively.
\begin{eqnarray}
&& R_c = \min_{\Delta_{P_i} \in \mathbb R^{N \times N},~ \Delta_b \in \mathbb R^{N \times M} }
\left\{\Vert \left[\Delta_{P_\dd} \cdots \Delta_{P_0}~~\Delta_b \right] \Vert_F ~|~ \right.  \nonumber \\
\label{rrc}
 && \left. (P_\dd+\Delta_{P_\dd},\dots P_0+\Delta_{P_0},b+\Delta_b) \textrm{ is uncontrollable} \right\}
\end{eqnarray}
Before we proceed further, there are several remarks in order.
\begin{remark}
{\rm
The perturbations $\Delta P_i \in \mathbb R^{N \times N},~ \Delta b \in \mathbb R^{N \times M}$
are allowed to be only real and hence the name \emph{real} radius of controllability. The case where complex
perturbations are allowed is not considered in this paper, though the theory proposed in this paper can be extended
very easily to accommodate for the complex case.
}
\end{remark}
\begin{remark}
{\rm
The norm of the perturbation matrices considered here is the Frobenius norm as opposed to the spectral norm as in previous works.
The advantage of Frobenius norm is that one can quantify the change in individual entries of $(P_\dd,\dots,P_1,P_0,b)$
 which has great physical advantage.
}
\end{remark}
\begin{remark}
\label{rmk_src}
{\rm
From the definition of real radius of controllability, it is clear that there is no restriction on the structure of the perturbations.
In many cases, such structure restriction on the perturbation matrices can arise naturally due to structure constraints of the systems.
In such a case, when the perturbations $\Delta P_i ~(i=0,1,\dots,\dd),~ \Delta b $ are structured perturbations, the radius of controllability,
for obvious reasons, is called as structured real radius of controllability. We note here, that the theory that we develop to
compute radius of controllability can be used to obtain the structured radius of controllability as well by slightly modifying
certain equations. This will be illustrated by several numerical examples. A particular case of interest that perturbation
to the leading coefficient matrix $P_\dd$ is not permitted is explored in the numerical examples presented later.
}
\end{remark}

In this paper, we transform the higher order system into a descriptor system, then we compute the radius of controllability
for the descriptor system only. The radius of controllability 
for descriptor system in \eqref{des}
  can be defined in similar way.
  \begin{eqnarray}
&& R_c = \min_{\Delta_E,\Delta_A \in \mathbb R^{n \times n},~ \Delta_B \in \mathbb R^{n \times m} }
\left\{\Vert \left[\Delta_E~~\Delta_A~~\Delta_B \right] \Vert_F ~|~ \right.  \nonumber \\
\label{rrc}
 && \left. (E+\Delta_E,A+\Delta_A,B+\Delta_B) \textrm{ is uncontrollable} \right\}.
\end{eqnarray}

The remainder of the paper is organized as follows: in the next section, we show the equivalence of computing the radius of controllability
with the computation of the nearest SLRA of a certain Toeplitz matrix.
In Section \ref{4}, we discuss an efficient algorithm (Structured Total Least Norm) to solve the SLRA problem. We illustrate the advantages of our
approach by demonstrating several numerical case studies in
Section \ref{5} before we conclude in Section \ref{6}.
\section{Mathematical Formulation}
\label{3}
In this section, we reformulate the problem of computing the radius of controllability of a linear time invariant
higher order system as the problem of computing the nearest SLRA of a certain Toeplitz matrix.
We first transform the higher order LTI system \eqref{higher} into its canonical form.
 The canonical form of equation \eqref{higher} can be written as a descriptor system
 \begin{equation}
 \label{des1}
 E\overset{.}{z}(t)=Az(t)+Bu(t),
\end{equation}
 where
$E=\begin{bmatrix}
P_\dd & 0 &\dots & 0 \\
0 & I_{N} &  \dots & 0 \\
&  & \ddots & \\
0 &  0 & \dots & I_{N}
\end{bmatrix} \in \mathbb R^{N\dd \times N\dd}$, \\
$A=\begin{bmatrix}
-P_{\dd-1} & -P_{\dd-2}  &\dots &  -P_1 & -P_0 \\
I_{N} & 0  & \dots &  0 & 0\\
0 & I_N  & \dots & 0 & 0\\
& & \ddots &  &\\
0  & 0 &  \dots & I_N & 0
\end{bmatrix} \in \mathbb R^{N\dd \times N\dd}$,
$B=\begin{bmatrix}
b^T &
0^T &
\cdots &
0^T
\end{bmatrix}\in \mathbb R^{N\dd \times M}$

Therefore System \eqref{des1} is same as System \eqref{des} with $n=N\dd$ and $m=M$.
Now we will state and prove a theorem which connects the controllability of the higher order system in \eqref{higher}
 and its canonical form in \eqref{des1}.
 Before proving the theorem, we state the following necessary lemma.
 \begin{lemma}\cite[~6.3-20]{kailath1980linear}
 \label{lemma1}
 Consider the $\tt r^{\rm{th}}$ order polynomial matrix $R(s)=R_{\tt r} s^{\tt r} +\dots+R_1 s+R_0$, where $R_i\in \mathbb R^{r_1 \times r_2}$ for $i=0,1,\dots,\tt r$. Then
 \begin{equation*}
 \mathcal{R}(s)=U(s)[sE_r-A_r]V(s)
 \end{equation*}
 where $\mathcal{R}(s)=\text{block diag}\{\underbrace{I_{r_2},\dots,I_{r_2},R(s)}_{\tt r\text{-terms}}\},$\\
 $E_r=\begin{bmatrix}
R_{\tt r} & 0 &\dots & 0 \\
0 & I_{r_2} &  \dots & 0 \\
&  & \ddots & \\
0 &  0 & \dots & I_{r_2}
\end{bmatrix}$,\\
$A_r=\begin{bmatrix}
-R_{\tt r-1} & -R_{\tt r-2}  &\dots &  -R_1 & -R_0 \\
I_{r_2} & 0  & \dots &  0 & 0\\
0 & I_{r_2}  & \dots & 0 & 0\\
& & \ddots &  &\\
0  & 0 &  \dots & I_{r_2} & 0
\end{bmatrix}$
  and $U(s),V(s)$ are \emph{unimodular}\footnote{A matrix $U(s)$ is called \emph{unimodular} if Det$[U(s)]$ is a non-zero constant.} matrices.
 \end{lemma}
\begin{t1}
\label{th4}
 The higher order system in \eqref{higher} is $C^\dd$-controllable if and only if
 its canonical form in \eqref{des1} is C-controllable.
\end{t1}
\emph{Proof:}
It is observed that if we can prove the equivalence of the conditions of Theorem \ref{th2} and Theorem \ref{th3}, then the proof is done.
 Now look at the structures of the matrices $E,A$ and $B$ shown in \eqref{des1}. From the structure of $[E~B]$,
 it is observed that $[E ~B]$ has full rank
 if and only if $[P_\dd~b]$ has full rank. Therefore condition (2) of Theorem \ref{th2} and Theorem \ref{th3} are equivalent.
 Now it is sufficient to prove only
that rank$([P(s)~b])=N$ if and only if rank$([sE-A~B])=n$.

Let $Q(s)=[P(s)~b]\in \mathbb R^{N \times (N+M)}$. It can be written that $Q(s)=Q_\dd s^\dd +\dots+Q_1 s+Q_0$,
where $Q_0=[P_0 ~b], ~Q_i=[P_i ~0]~\text{for}~i=1,2,\dots,\dd.$ Now from Lemma \ref{lemma1}, we get there exist unimodular matrices
$U(s)$ and $V(s)$ such that $\mathcal{Q}(s)=U(s)[sE_1-A_1]V(s)$
 where $\mathcal{Q}(s)=\text{block diag}\{\underbrace{I_{(N+M)},\dots,I_{(N+M)},Q(s)}_{\dd\text{-terms}}\},$\\
$E_1=\textrm{diag}\begin{bmatrix}
Q_\dd & & I_{(N+M)} & \cdots & I_{(N+M)}
\end{bmatrix}$, and\\
$A_1=\begin{bmatrix}
-Q_{\dd-1} & -Q_{\dd-2}  &\dots &  -Q_1 & -Q_0 \\
I_{(N+M)} & 0  & \dots &  0 & 0\\
0 & I_{(N+M)}  & \dots & 0 & 0\\
& & \ddots &  &\\
0  & 0 &  \dots & I_{(N+M)} & 0
\end{bmatrix}$\\
Therefore $\mathcal{Q}(s)$ has full rank if and only if $[sE_1-A_1]$ has full rank for all finite $s \in \mathbb C$.
Now substitute $Q(s)=[P(s)~b]$ in $\mathcal{Q}(s)$, then after some elementary row/column operations we get\\
\[ \mathcal{Q}(s) \equiv \left[\begin{array}{cccc|cccc}
I_N & 0 &\dots & 0 & 0  & 0 & \dots & 0\\
0 & I_{N} &  \dots & 0 & 0 & 0 & \dots &  0\\
&  & \ddots & & &  & &\\
0 &  0 & \dots & P(s) & 0 & 0 & \dots & b\\
\hline
0 & 0 &\dots & 0 & I_M & 0 & \dots & 0\\
0 & 0 &  \dots & 0 & 0 & I_M & \dots &  0\\
&  & \ddots & & &  & \ddots &\\
0 &  0 & \dots & 0 & 0 & 0 & \dots & I_{M}
\end{array}
\right] \] \\
Now from this structure it is clear that $\mathcal{Q}(s)$ has full rank if and only if $[P(s)~b]$ has full rank for all finite $s \in \mathbb C$.
Again substitute $Q_0=[P_0 ~b], ~Q_i=[P_i ~0]~\text{for}~i=1,2,\dots,\dd$ in $sE_1-A_1$. Similarly after some elementary
row/column operations we get
$ sE_1-A_1 \equiv$  \[
\left[\begin{array}{c|ccccc}
[sE-A~B]  & 0  & 0 & 0 & \dots & 0\\
\hline
0 &  -I_M & sI_M & 0 & \dots & 0\\
0 &  0 & -I_M & sI_M & \dots &  0\\
& & &  & \ddots &\\
0 & 0 & 0 & \dots & -I_M & sI_{M}
\end{array}
\right] \] \\
From the structure of the matrix, it is also clear that $sE_1-A_1$ has full rank if and only if $[sE-A~B]$ has full rank for all finite
$s \in \mathbb C.$ Therefore combining the results we get for all finite $s \in \mathbb C$,
\begin{align*}
[P(s)~b]\text{ has full rank } \Leftrightarrow & \mathcal{Q}(s)\text{ has full rank} \\
\Leftrightarrow & [sE_1-A_1] \text{ has full rank} \\
\Leftrightarrow & [sE-A~B] \text{ has full rank.}
\end{align*}
Hence the result is proved.
\QEDB
\begin{remark}
\label{rmk:highertolower}
{\rm
From Theorem \ref{th4}, it directly follows that nearest uncontrollable system to the higher order system \eqref{higher} can
be obtained from the nearest uncontrollable system to the descriptor system \eqref{des1}. Suppose $(\tilde{E},\tilde{A},\tilde{B})$
 is the nearest uncontrollable system to System \eqref{des1} satisfying the additional constraint that $(\tilde{E},\tilde{A},\tilde{B})$ and
 $(E,A,B)$ have the same structures. Therefore from \eqref{des1} it can be observed that the perturbations are allowed only in the blocks
 $(P_\dd,\dots,P_1,P_0,b)$ in the perturbed $(E,A,B)$. Suppose the perturbed matrices $(\tilde{E},\tilde{A},\tilde{B})$ follows the
 perturbed matrices $(\tilde{P}_\dd,\dots,\tilde{P}_1,\tilde{P}_0,\tilde{b})$. Then we can conclude that the corresponding
 perturbed higher order system $(\tilde{P}_\dd,\dots,\tilde{P}_1,\tilde{P}_0,\tilde{b})$ is the nearest uncontrollable system to \eqref{higher}. Therefore, in order to compute the real radius of controllability of higher order system, it is sufficient to compute the structured real radius of controllability (the radius of controllability obtained by respecting the structure of perturbation in the parameter matrices) of the equivalent descriptor system.
}
\end{remark}

In light of condition (2) of Theorem \ref{th3},
the problem of computing the real radius of controllability of System \eqref{des} can be split into
two following subproblems:
\begin{problem}
\label{p1}
\begin{align*}
r_1=&\displaystyle \min_{\Delta_E,~\Delta_A,~ \Delta_B}  \left\{\| [\Delta_E ~~\Delta_A ~~\Delta_B] \|_F, \text{ such that } \right.\\
&\left.\text{rank}[s(E+\Delta_E)-(A+\Delta_A)~~B+\Delta_B]<n \right. \\
&\left.\text{ for some }s \in \mathbb C \right\}
\end{align*}
\end{problem}
\begin{problem}
\label{p2}
\begin{align*}
&r_2=\displaystyle \min_{\Delta_E,~ \Delta_B}  \left\{\| [\Delta_E ~~\Delta_B] \|_F, \right.\\
&\left.\text{ such that, } \text{rank}[(E+\Delta_E) ~~(B+\Delta_B)]<n \right\}
\end{align*}
\end{problem}
Therefore it can be observed that radius of controllability of the system \eqref{des} can be found by choosing
$R_c=\min \{r_1,r_2\}$ with the understanding that
the system will be uncontrollable if any of the matrices among $[sE-A~~ B]$ and $[E ~~B]$ loses its rank.
In \cite{simon}, these two problems were solved individually to find $r_1$ and $r_2$.
We now present a result which states that controllability of System \eqref{des} can be shown to be equivalent to the full rank
property of a certain Toeplitz structured matrix.
\begin{t1}
\label{th5}
\emph{\cite[Theorem 2-2.1]{dai}}
Consider the descriptor system given in \eqref{des}.
Then the following statements are equivalent.
\begin{enumerate}
\item The system in \eqref{des} is C-controllable.
\item
\label{equivalenceCR}
$ \mathcal{C}(E,A,B) \in \mathbb R^{n^2 \times n(n+m-1)}$ is full row-rank matrix, where,
\begin{eqnarray}
&& \mathcal{C}(E,A,B) \nonumber \\
&=&
\label{Toeplitz}
\begin{bmatrix}
-A & B &    &   &        &    &  & \\
 E & 0 & -A & B &        &    &  & \\
   &   & E  & 0 &        &    &  & \\
   &   &    &   &\ddots  &    &  & \\
   &   &    &   &        & -A & B & \\
   &   &    &   &        & E  & 0 & B
\end{bmatrix}
\end{eqnarray}
\end{enumerate}
\end{t1}

The equivalence of C-controllability in Statement \ref{equivalenceCR} of Theorem \ref{th5} is the key to formulate the
computation of the real radius of controllability as a structured low rank approximation problem. From Theorem \ref{th5}, we immediately conclude the following theorem which is the main result of the paper.
\begin{t1}
 \label{slra-rrc}
 Let an LTI descriptor system be given as in \eqref{des}. Then computing the radius of controllability, denoted by $R_c$, of the system
 in \eqref{des}, is equivalent to computing nearest SLRA $~\mathcal{C}(\hat{E},\hat{A},\hat{B})$ to $\mathcal{C}(E,A,B)$, introduced in Statement \ref{equivalenceCR} of Theorem \ref{th5}.
 In fact,
 \begin{equation}
  \label{rrc_formula}
  R_c = \left \Vert \left[ \Delta E ~~ \Delta A ~~ \Delta B \right] \right \Vert_F.
 \end{equation}
 where, $\hat{E}=E+\Delta E, ~\hat{A}=A+\Delta A,~ \hat{B}=B+\Delta B$.
\end{t1}

\emph{Proof:}
Suppose that the given system $(E,A,B)$
is C-controllable, then the statement \ref{equivalenceCR} states that the Toeplitz structured matrix
$\mathcal{C}(E,A,B) \in \mathbb R^{n^2 \times n(n+m-1)}$ is a full row rank matrix. In order to
compute the real radius of controllability, we seek the perturbations $\Delta E, \Delta A, \Delta B$ of appropriate dimensions such that
the perturbed system $(E+\Delta E, A+\Delta A, B+\Delta B)$ is uncontrollable. Then Statement \ref{equivalenceCR} in Theorem \ref{th5}
would imply that the Toeplitz structured matrix associated with the perturbed system,
denoted as $\mathcal{C} (E+\Delta E, A+\Delta A, B+\Delta B)$ as in \eqref{perturbedT}, is not a full row rank matrix.
\begin{eqnarray}
&& \mathcal{C} (E+\Delta E, A+\Delta A, B+\Delta B) \nonumber \\
\label{perturbedT}
&&=
\begin{bmatrix}
-\hat{A} & \hat{B} &    &   &        &    &  &\\
 \hat{E} & 0 & -\hat{A} & \hat{B} &        &    &   &\\
   &   & \hat{E}  & 0 &        &    &  &\\
   &   &    &   &\ddots  &    & &\\
   &   &    &   &        & -\hat{A} & \hat{B}  &\\
   &   &    &   &        & \hat{E}  & 0  & \hat{B}
\end{bmatrix}
\end{eqnarray}
where $\hat{A} = A +\Delta A$, $\hat{E} = E + \Delta E$ and $\hat{B} = B + \Delta B$. Clearly, $\mathcal{C} (E+\Delta E, A+\Delta A, B+\Delta B)$
is a structured perturbation of $\mathcal{C} (E, A, B)$.

Conversely, given a system $(E,A,B)$, construct the Toeplitz structured matrix as in \eqref{Toeplitz}. If we are able to compute some perturbations
$(\Delta E, \Delta A, \Delta B)$ of appropriate sizes such that the perturbed Toeplitz matrix in \eqref{perturbedT} is not full row rank, then
the corresponding system $(E+\Delta E, A+\Delta A, B+\Delta B)$ is clearly uncontrollable. Let these perturbations, in particular, have the following
property.
\begin{eqnarray}
 && \left \Vert \left[ \Delta E~~\Delta A~~ \Delta B \right] \right \Vert_F
  = \min_{\Delta_E, \Delta_A, \Delta_B}
  \left\{ \left \Vert \left[ \Delta_E~~\Delta_A~~ \Delta_B \right] \right \Vert_F ~|~ \right. \nonumber \\
 \label{minperturb}
&& \left. \mathcal{C} (E+\Delta_E, A+\Delta_A, B+\Delta_B) \textrm{ is not full row rank}\right\}
\end{eqnarray}
In such a case, the matrix $\mathcal{C} (E+\Delta E, A+\Delta A, B+\Delta B)$ is the nearest structured low rank approximation to the matrix
$\mathcal{C} (E, A, B)$. Clearly, from \eqref{minperturb} both these matrices have the same structure. Further, not only the matrix
$\mathcal{C} (E+\Delta E, A+\Delta A, B+\Delta B)$ is of lower rank than $\mathcal{C} (E, A, B)$, but is also the nearest with the low rank property.
This justifies why this perturbed matrix is called as the nearest Structured Low Rank Approximation (SLRA) of the given matrix.
Thus, from \eqref{rrc} and \eqref{minperturb}, it is clear that the problem of computing the real radius of controllability of the system
in \eqref{des} is equivalent to computation of the nearest SLRA $~\mathcal{C}(\hat{E},\hat{A},\hat{B})$ of $\mathcal{C}(E,A,B)$. Also clearly $R_c = \left \Vert \left[ \Delta E ~~ \Delta A ~~ \Delta B \right] \right \Vert_F,$ where, $\hat{E}=E+\Delta E, ~\hat{A}=A+\Delta A,~ \hat{B}=B+\Delta B$. \QEDB

%
\begin{remark}
{\rm
Note that we are using Statement \ref{equivalenceCR} in Theorem \ref{th5} to compute the radius of controllability. This enables us to compute
the radius of controllability directly as opposed to the other approaches in the literature where $r_1$ and $r_2$ are computed as in Problem \ref{p1}
and Problem \ref{p2} and the minimum of these is declared to be $R_c$.
}
\end{remark}
\begin{remark}
{\rm
From the definition and Theorem \ref{slra-rrc}, we can understand the radius of controllability in another way. Theorem \ref{slra-rrc}, in fact,
says that there a ball around the system $(E,A,B)$ with radius equal to $R_c$ such that every system inside this ball is C-controllable. The larger
radius of this ball ensures higher robustness with respect to controllability.
On the other hand a small value of $R_c$ indicates that an uncontrollable
system is nearby indicating the system may be susceptible to small perturbations.
}
\end{remark}

In the literature, the problem of computing the nearest SLRA to a given matrix is studied extensively.
We will explain the problem of SLRA formally and a procedure to solve this
problem in the next section.

As pointed out in Remark \ref{rmk_src} and Remark \ref{rmk:highertolower}, it may so happen that certain physical constraints may impose structural constraints on the perturbations
while computing the nearest uncontrollable system to \eqref{des}. We now define \emph{Structured Real Radius of Controllability} for descriptor
systems under the assumption that these structure constraints on the perturbations are linear or affine. Let $\mathcal{S}_E \subseteq \mathbb{R}^{n \times n}$
be a linear subspace (or affine subset) which characterizes the perturbations allowed in the system matrix $E$. Similarly,
let $\mathcal{S}_A \subseteq \mathbb{R}^{n \times n}$ and $\mathcal{S}_B \subseteq \mathbb{R}^{n \times m}$ be the linear subspaces
(or affine subset) which characterize
the perturbations to be given to the system matrices $A$ and $B$ respectively. Then the Structured Real Radius of Controllability is
defined as follows.
\begin{eqnarray}
&& SR_c = \min_{\Delta_E \in \mathcal{S}_E,~\Delta_A \in \mathcal{S}_A, ~ \Delta_B \in \mathcal{S}_B }
\left\{\Vert \left[\Delta_E~~\Delta_A~~\Delta_B \right] \Vert_F ~|~ \right.  \nonumber \\
\label{Srrc}
 && \left. (E+\Delta_E,A+\Delta_A,B+\Delta_B) \textrm{ is uncontrollable} \right\}.
\end{eqnarray}
It should be clear that by performing the same calculations as in the case of real radius of controllability,
we can show that the problem of computing the structured real radius of controllability can be posed as an SLRA problem. For the sake
of completeness we state the following theorem.
\begin{t1}
 \label{slra-Srrc}
 Let an LTI descriptor system be given as in \eqref{des}. Let $SR_c$ denote the structured real radius of controllability
 of the system in \eqref{des}. Then $SR_c$ can be computed by computing the nearest SLRA
 $\mathcal{C} (E+\Delta E, A+\Delta A, B+\Delta B)$ to $\mathcal{C}(E,A,B)$ by solving the following optimization problem.
 \begin{eqnarray*}
 \label{srrc_formula}
 && SR_c =  \min_{\Delta_E \in \mathcal{S}_E, \Delta_A \in \mathcal{S}_A, \Delta_B \in \mathcal{S}_B }
  \left\{ \left \Vert \left[ \Delta_E~~\Delta_A~~ \Delta_B \right] \right \Vert_F ~|~ \right. \nonumber \\
&& \left. \mathcal{C} (E+\Delta_E, A+\Delta_A, B+\Delta_B) \textrm{ is not full row rank}\right\}.
\end{eqnarray*}
\end{t1}
\emph{Proof:} In the Theorem \ref{slra-rrc}, the nearest SLRA was computed to $\mathcal{C}(E,A,B)$ for any $E,A \in \mathbb R^{n \times n},~B \in \mathbb R^{n \times m}$, but here the same thing has to be done for the structured matrices $E \in \mathcal{S}_E,~A \in \mathcal{S}_A,~B \in \mathcal{S}_B$. Therefore, the proof is similar as the proof of Theorem \ref{slra-rrc}, we omit this. \QEDB

We now give the algorithm to compute the radius of controllability for a given higher order system.
\begin{algorithm} \emph{Structured real radius of controllability}\\
\underline{\emph{Input:}} Higher order system $(P_\dd,\dots,P_1,P_0,b)$\\
\underline{\emph{Output:}} The structured real radius of controllability, $SR_c$\\
\emph{Step 1:} Construct structured matrices $(E,A,B)$ given in Eq.\eqref{des1} \\
\emph{Step 2:} Construct the Toeplitz structured matrix $\mathcal{C}(E,A,B)$ from the triple $(E,A,B)$ as in \eqref{Toeplitz}.\\
\emph{Step 3:} Compute nearest SLRA $\mathcal{C}(E+\Delta E,A +\Delta A,B + \Delta B)$ to $\mathcal{C}(E,A,B)$
as in \eqref{minperturb} (See Algorithm 4.1). \\
\emph{Step 3:} Compute $SR_c = \left \Vert \left[ \Delta E ~~ \Delta A ~~ \Delta B \right] \right \Vert_F$.
\end{algorithm}

\section{Structured Low Rank Approximation (SLRA)}
\label{4}
In this section, we state the problem of computing the nearest SLRA to a given matrix and explain the Structured Total Least
Norm (STLN) algorithm to compute the nearest SLRA to the given matrix.
Let $X_* \in \mathbb R^{M \times N}$ be a given matrix with a particular affine structure.
We want to compute the nearest low rank approximation $X$ which inherits the same structure as $X_*$.
Therefore, the SLRA problem is an optimization problem where the distance between $X_*$ and $X$
has to be minimized. Mathematically, the SLRA problem can be formulated as
\begin{equation}
\label{slra}
\begin{aligned}
\displaystyle \min_{X \in  \mathcal{X} \cap \Omega } {\|X-X_*\|_F}, \\
\end{aligned}
\end{equation}
\text{where,}\\
$\begin{aligned}
&\Omega \in \mathbb R^{M \times N} \text{ is a space (or set) of particular linear (or affine)} \\
&~~~~~~~~~~~~~~ \text{ structured matrices},\\
&\mathcal{X} := \left \{X \in \mathbb R^{M \times N} ~~|~~ \text{rank}(X)  <  \text{rank}(X_*) \right \}.
\end{aligned}$
In literature, many techniques are reported to solve SLRA problem (see for instance \cite{s7,s4,s10} and the references therein).
In this paper, we will use Structured Total Least Norm (STLN) algorithm to compute the nearest SLRA \cite{s7}.
\subsection{Structured total least norm (STLN)}
Suppose the given matrix $X_*$ is partitioned like $X_*=[Y ~ ~ y]$, where $y$ is the last column of $X_*$.
The target is to solve the system $Y z\approx y$ approximately.  In other words, we aim to find the minimum perturbations in
$Y$ and $y$ in such a way that $(Y+E_1)z =(y+f_1),$ where $[E_1~~ f_1]$ belongs to the linear structure of the space $\Omega$.
Let $\{\Phi_1, \Phi_2, \ldots \Phi_{\ell}\}$ be an ordered basis of the linear space of $\Omega$ and $\alpha \in \mathbb{R}^{\ell}$ be the
\emph{co-ordinate} vector of $X_*$ with respect to this basis. Then the STLN problem can be stated as:
\begin{equation}
\displaystyle \min_{\alpha,z} \left \Vert \begin{bmatrix}
\omega r(\alpha, z) \\
D \alpha
\end{bmatrix} \right \Vert_2
\end{equation}
where $r(\alpha, z)=y-(Y+E_1)z$ is the residual of the system, $\omega$ is any large weight and $D$ is a suitably chosen weight matrix.
In our case, we assume $D=I$. Let a small change in $\alpha$ be denoted as $\Delta \alpha$ and a small change in $z$ be $\Delta z$,
and a small perturbation in the error matrix $E_1$ be $\Delta E_1$. Then the perturbed residual becomes
\begin{align*}
r(\alpha +\Delta & \alpha, z+\Delta z)=(y+f_1)-(Y+E_1+\Delta E_1)(z+\Delta z)\\
& \approx (y+f_1)-(Y+E_1)z-(Y+E_1)\Delta z-\Delta E_1 z \\
& \approx r(\alpha, z)+(f_1-\Delta E_1 z)-(Y+E_1)\Delta z
\end{align*}
Compute a matrix $P$ such that $f_1=P\Delta \alpha$, and also compute a structured matrix $S$ dependent on $z$ which satisfies $S \Delta \alpha =\Delta E_1 z$. Formation of the matrix $S$ is explained in \cite{rosen}. Then $r(\alpha +\Delta \alpha, z+\Delta z)$ becomes $r(\alpha, z)-(S-P)\Delta \alpha-(A_1+E_1)\Delta z$. In this notation, the STLN problem is reformulated as
\begin{equation}
\label{stln}
\displaystyle \min_{\Delta \alpha, \Delta z} \left \Vert \begin{bmatrix}
\omega (S-P) & Y+E_1 \\
I & 0
\end{bmatrix}
\begin{bmatrix}
\Delta \alpha \\
\Delta z
\end{bmatrix}+
\begin{bmatrix}
-\omega r \\
\alpha
\end{bmatrix} \right \Vert_2
\end{equation}
Therefore the Algorithm for STLN can be written as:
\begin{algorithm} \emph{STLN}\\
\label{STLN_Algo}
\underline{\emph{Input}}: Matrices $Y,~ y,~ \text{error tolerance } \epsilon >0$\\
\underline{\emph{Output}}: Error matrix $[E_1 ~~f_1]$, with least norm.\\
Step 1: Set $\omega$ to be a large number, $E_1=0, f_1=0$ and find $z$ from $\min_z \|y-Yz\|$, and construct $S$ from $z$.\\
Step 2: Set $r=y-Yz$\\
Step 3: Repeat \\
$~~~~~~~~~~$       a. Solve the minimization problem in \eqref{stln}.\\
$~~~~~~~~~~$       b. Set $z:=z+\Delta z, \alpha :=\alpha + \Delta \alpha $\\
$~~~~~~~~~~$       c. Construct $[E_1~~ f_1]$ from $\alpha$ and $S$ from $z$ \\
$~~~~~~~~~~$       d. $r=(y+f_1)-(Y+E_1)z$ \\
until ($\| \alpha \|, \| z \| < \epsilon$)
\end{algorithm}
\begin{remark}
{\rm
In every iteration of Algorithm \ref{STLN_Algo}, the optimization problem in \eqref{stln} is solved.
Note that, this optimization problem is a least square problem. There are many efficient methods like $QR$ decomposition
to solve the least square problem efficiently. Further, it is clear from  \eqref{stln} that matrices involved have a
particular structure of zeros which can be exploited to solve this least square problem much more efficiently.
}
\end{remark}
\begin{remark}
{\rm
It should be noticed that the optimization problem in SLRA formulation has many local solutions (as explained in \cite{gillard2013optimization}).
STLN algorithm does not guarantee the global solution to this optimization problem.
Also, STLN may be very sensitive with respect to the partition of the matrix we are considering at the
beginning of the formulation. In the formulation, we have partitioned the matrix by its last column as $X_*=[Y ~ ~ y]$,
where $y$ is the last column. However, instead of chosing the last column, any other column of the matrix $X_*$ can be chosen at
the time of partition. Though theoretically chosing any partition of the matrix $X_*$ is equivalent and should lead to the same result, in practice
it is observed that this affects the numerical results. Therefore care must be exercised while choosing the partition.
Another important factor in STLN is the choice of the constant $\omega$. This choice may also affect the rate of convergence as well the
local minimum to which STLN converges.
}
\end{remark}
\begin{remark}
{\rm
As stated earlier, we now explain a way in which we can compute the structured radius of controllability.
Observe that the structure constraint imposed on the perturbation matrices can be characterized by a linear subspace.
In the STLN formulation, we then chose a basis for this linear structure space accordingly.
Once this is done, then the rest of the calculation remains same. This shows that, in this formulation, any linear structure
restriction on perturbations can be very easily incorporated. This is a great advantage of the proposed formulation.
}
\end{remark}
\section{Numerical Examples}
\label{5}
In this section, we study some numerical examples to compute the structured radius of controllability of higher order, descriptor systems.
\begin{Example}
This is an example taken from the paper \cite{simon}. The descriptor system is given by \begin{equation}
\label{ex1}
E\overset{.}{x}(t) = Ax(t) + Bu(t)
\end{equation}
 where, the system matrices are given by \\ $E=\begin{bmatrix}
1.8 & 0 & 0 \\
0 & 0.34 & 0 \\
0 & 0 & 0\\
\end{bmatrix},
A=\begin{bmatrix}
2 & -0.91 & -0.088 \\
0.19 & 0.25 & 0.51 \\
0.64 & 0.31 & -0.59
\end{bmatrix},\\
B=\begin{bmatrix}
-0.63 \\
0.53 \\
-0.58
\end{bmatrix}.$
The system \eqref{ex1} is controllable, our interest is to compute the nearest uncontrollable system.
In this case it is assumed that no perturbation in the matrix $E$ is allowed i.e., $\Delta E=0$.
According to the theory, STLN technique (seting, $\omega=10^{13}$) is used to compute the SLRA of the matrix
$\mathcal{C}(E,A,B) \in \mathbb R^{9 \times 9}$, constructed from the given $(E,A,B)$.  In this case, we get
the structured radius of controllability $SR_c=0.3436$, (we have taken 2-norm of the perturbation matrix in order to facilitate
the comparison with the results in \cite{simon})
and the corresponding nearest uncontrollable system is $\hat{E}\overset{.}{x}(t) = \hat{A}x(t) + \hat{B}u(t)$, where \\
$\hat{A}=A+\Delta A=\begin{bmatrix}
 2.0000 & -0.9100 & -0.3333 \\
   0.1900  & 0.2500 & 0.4362 \\
   0.6400  & 0.3100 & -0.4887
\end{bmatrix}\\
\hat{B}=B+\Delta B=\begin{bmatrix}
-0.4471 \\
   0.5850 \\
  -0.6555
 \end{bmatrix}$. It can be observed that the result exactly matches with the result shown in \cite{simon}.
\end{Example}
\begin{Example}
A parametric family of descriptor systems is taken as $E\overset{.}{x} = Ax + Bu$, where the system matrices are given as \\
$E=\begin{bmatrix}
0 &2.1 &0 \\
1 &0 &0\\
0 &0 &0
\end{bmatrix}
A=\begin{bmatrix}
1 & 3 & 0\\
2 & 1 & 1\\
3 & 1 & 5
\end{bmatrix}
B=\begin{bmatrix}
1\\
0\\
\delta
\end{bmatrix}, \delta \in \mathbb R$\\
We have computed the radius of controllability for different values of $\delta$ by our approach, considering no
perturbation in $E$. According to the structure of the matrix $[E~~ B]=\begin{bmatrix}
0 &2.1 &0 &1\\
1 &0 &0 &0\\
0 &0 &0 & \delta
\end{bmatrix}$, it is clear that rank$[E~~B]$ is determined by the parameter $\delta$. When the value of $\delta$ is
very close to zero, the matrix $[E~~B]$ is going to loose its rank, equivalently the system is going
to become uncontrollable. Therefore according to our approach, using the matrices $E,A,B$  we constructed
the matrix $\mathcal{C}(E,A,B) \in \mathbb R^{9 \times 9}$, and computed the SLRA using STLN technique. The final results
are shown in the following table.
\begin{table}[htbp]
\centering
\begin{tabular}{|c||ccccccc|}
\hline
$\delta$ & 1& 0.6 & 0.4& 0.2& 0.1 & 0.01 & 0\\
\hline
$SR_c$ & 0.3193 & 0.3820 & 0.4132 & 0.2000 &  0.1 & 0.01 & 0 \\
\hline
\end{tabular}
\end{table}
\end{Example}
\begin{Example}
The following example, based on an electric circuit system, is adopted from \cite{dai}. Consider the circuit given in the Figure~\ref{fig:circuit}.
\begin{figure}[!h]
\centering
\includegraphics[height= 3cm,width=7.5cm]{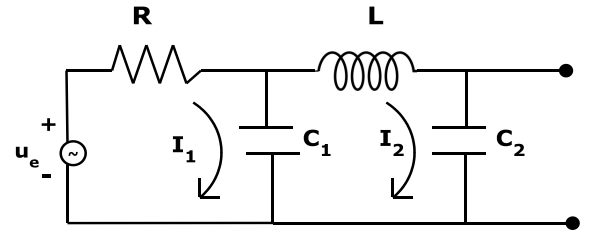}
\caption{In the electric circuit, $R$ stands for resistor, $L$ stands for inductor and $C_1, ~C_2$ stand for capacitors.}
\label{fig:circuit}
\end{figure}
Choose the state variable as $x=[u_{C_1}~~ u_{C_2}~~ I_2 ~~I_1]^T$, where $u_{C_1}, u_{C_2}, I_2, I_1$ are the voltages of capacitors:
$C_1,~C_2$ and the currents flowing over them. According to Kirchoff's laws, the following linear descriptor system can be obtained:
\begin{equation}
\label{ec}
\begin{bmatrix}
C_1 & 0 & 0 & 0\\
0 & C_2 & 0 & 0\\
0 & 0 & -L & 0\\
0 & 0 & 0 & 0
\end{bmatrix} \overset{.}{x}= \begin{bmatrix}
0 &0 &0 &1\\
0& 0 &1 &0\\
-1& 1 &0& 0\\
1 &0& 0& R
\end{bmatrix}x+\begin{bmatrix}
0\\
0\\
0\\
-1
\end{bmatrix}u_e
\end{equation}
For this system, we have,  rank$[sE-A~~B]=\text{rank}\begin{bmatrix}
sC_1 & 0 & 0 & -1 & 0\\
0 & sC_2 & -1 & 0 & 0\\
1 & -1 & -sL & 0 & 0\\
-1 & 0 & 0 & -R & -1
\end{bmatrix}=4$,  $\forall s \in \mathbb C $.
Also rank$[E~~B]$=rank$\begin{bmatrix}
C_1 & 0 & 0 & 0 & 0\\
0 & C_2 & 0 & 0 & 0\\
0 & 0 & -L & 0 & 0\\
0 & 0 & 0 & 0 & -1
\end{bmatrix}=4$ for any non-zero values of $C_1,C_2$ and $L$.
Therefore the system \eqref{ec} is controllable for any values of $R$ and for any non-zero values of $C_1,C_2$ and $L$.
It is clear from the structure of the matrices that $[sE-A~~B]$ can not lose rank for any values of $C_1,C_2,L,R$,
however the only way for the matrix $[E~~B]$ to lose its rank is by setting any of the entries $C_1,C_2,L$ to zero.
Therefore, the nearest uncontrollable system can be obtained by setting the smallest entry of $C_1,C_2,L$ to zero.
The radius of controllability is then $\min \{C_1,C_2,L\}$.

We can verify the result by computing the structured radius of controllability $SR_c$ in the approach we proposed here.
Using the matrices $E,A,B$ we obtain $\mathcal{C}(E,A,B) \in \mathbb R^{16 \times 16}$, whose nearest SLRA needs to be computed.
However it should be remembered that the perturbations are permitted in $C_1,C_2,L$ and $R$ only, the other elements should remain
unchanged. Therefore, this is the case of finding structured radius of controllability, where not all the entries system matrices
are permitted to be perturbed.
For different values of $C_1,C_2,L$ and $R$, our results are given in Table \ref{tab:tab3}, from which the accuracy and efficiency of
our technique is proved.
\begin{table}[htbp]
\caption{}
\centering
\begin{tabular}{|c||c|c|}
\hline
$(C_1,C_2,L,R)$ & $\min \{C_1,C_2,L\}$ & $SR_c$ \\
\hline
(1,~1,~1,~1) & 1 & 0.9997 \\
(2,~1.5,~3,~1) & 1.5 & 1.4998 \\
(2,~3.5,~1.2,~4) & 1.2 & 1.2000 \\
(0.0001,~0.1,~10,~3) & 0.0001 & 0.0001 \\
(8,~0.01,~0.1,~4) & 0.01 & 0.0100 \\
\hline
\end{tabular}
\label{tab:tab3}
\end{table}
\end{Example}
\begin{Example}
In this example, we have considered ten randomly generated linear descriptor systems $(E,A,B)$ of different dimensions. In each of these
cases, we have computed the structred real radius of controllability where we do not allow perturbations in $E$.
The purpose of this example is to show the computational efforts
required to compute the radius of controllability in terms of time and number of iterations required in STLN algorithm.
We note the average number of iterations in STLN (taking $\omega=10^{8}$, $\epsilon=10^{-3}$) and average time required to compute the nearest SLRA.
These results are summarized in Table \ref{tab:tab4}.
\begin{table}[htbp]
\caption{}
\centering
\begin{tabular}{|c||c|c|}
\hline
$(n,m)$ & Avg. Time (Seconds) & Avg. No. of iterations \\
\hline
(5,~1) & 0.0166 & 8.3 \\
(7,~2) & 0.0957 & 18.5 \\
(10,~3) & 0.6242 & 23.9 \\
(15,~4) & 2.4168 & 22.6 \\
(20,~5) & 11.1736 & 20.5 \\
\hline
\end{tabular}
\label{tab:tab4}
\end{table}
It should also be noted that in the case of multi-input descriptor systems, the matrix $C(E,A,B)$ is a wide matrix. Therefore, we use transpose
of the $C(E,A,B)$ and implement STLN as discussed.
\end{Example}
\begin{Example}
(\textbf{A quadratic brake model}) In \cite{simon,mengi2008estimation}, the vibrations of a drum brake system are studied. The model consists of a
rotating disc with mass $m$, radius $r$, and two contact points subject to friction and is described by the quadratic equation
\begin{equation}
\label{brake}
M\overset{..}x(t)+K(\mu)x(t)=f(t)
\end{equation}
where the mass matrix $M$ is given by
\begin{equation*}
\begin{bmatrix}
m & 0\\
0 & m
\end{bmatrix},
\end{equation*}
and the stiffness matrix $K(\mu)$ is
\begin{equation*}
k\begin{bmatrix}
(\sin\gamma +\mu\cos\gamma)\sin\gamma & -\mu-(\sin\gamma+\mu\cos\gamma)\cos\gamma \\
(\mu\sin\gamma -\cos\gamma)\sin\gamma &1+(\mu\sin\gamma+\cos\gamma)\cos\gamma
\end{bmatrix},
\end{equation*}
$\mu$ is the friction coefficient, and $k$ and $\gamma$ are other parameters specific to the model. As in \cite{simon,mengi2008estimation},
suppose the force on the brake system has just the vertical component determined by the input
\begin{equation}
f(t)=\begin{bmatrix}
0 \\1
\end{bmatrix}u(t)
\end{equation}
System \eqref{brake} is a second order system, which can be transformed into its first order canonical form
 \begin{equation}
\label{robotic1}
\begin{bmatrix}
M & 0\\
0 & I
\end{bmatrix}\overset{.}z=
\begin{bmatrix}
0 & -K(\mu) \\
I & 0
\end{bmatrix}\overset{.}z+
\begin{bmatrix}
B\\
0
\end{bmatrix}u,
\end{equation} where $B=\begin{bmatrix}0 \\ 1\end{bmatrix}$.
For the parameters $m=5,~k=1,~\gamma=\frac{\pi}{100}$, we computed structured radius of controllability for different values of $\mu$. It is assumed that
the mass matrix $M$ is not supposed to get perturbed, only perturbations are allowed in the matrices $K$ and $B$. The results are compared with the
corresponding results from \cite{simon} (see Table \ref{rotating_disc}). It is noticed that the results are exactly same for most of
the cases, which proves the efficiency of our algorithm.
\begin{table}[h]
\caption{}
\centering
\begin{tabular}{|c||c|c|}
\hline
$\mu$  & \multicolumn{2}{|c|}{$SR_c$} \\
\hline
&  Using SLRA & According to \cite{simon} \\
\hline
0.05 & 0.0587 & 0.0587\\
0.1 & 0.1031 & 0.1031\\
0.15 & 0.1470 & 0.1470\\
0.2 & 0.1901 & 0.1901\\
0.5 & 0.4227 & 0.4227\\
1 & 0.6813 & 0.6811\\
10 & 0.9959 & 0.9959\\
100 & 1.0000 & 1.0000\\
1000 & 1.0000 & 1.0000\\
\hline
\end{tabular}
\label{rotating_disc}
\end{table}
\end{Example}
\section{Concluding Remarks}
\label{6}
In this paper, we considered the problem of computing the
structured radius of controllability for higher order systems. It was shown that this problem is equivalent to the problem of
computing the nearest SLRA of a certain Toeplitz structured matrix obtained from the system matrices. It is
also shown that the structure constraint on the perturbations can be inherently modeled in the SLRA
formulation making the computation of structured radius of controllability on par with the radius of controllability.
Several numerical case studies were presented to highlight advantages of the proposed algorithm.
\bibliographystyle{IEEEtran}
\bibliography{rrc2}

\begin{thebibliography}{10}
\providecommand{\url}[1]{#1}
\csname url@samestyle\endcsname
\providecommand{\newblock}{\relax}
\providecommand{\bibinfo}[2]{#2}
\providecommand{\BIBentrySTDinterwordspacing}{\spaceskip=0pt\relax}
\providecommand{\BIBentryALTinterwordstretchfactor}{4}
\providecommand{\BIBentryALTinterwordspacing}{\spaceskip=\fontdimen2\font plus
\BIBentryALTinterwordstretchfactor\fontdimen3\font minus
  \fontdimen4\font\relax}
\providecommand{\BIBforeignlanguage}[2]{{%
\expandafter\ifx\csname l@#1\endcsname\relax
\typeout{** WARNING: IEEEtran.bst: No hyphenation pattern has been}%
\typeout{** loaded for the language `#1'. Using the pattern for}%
\typeout{** the default language instead.}%
\else
\language=\csname l@#1\endcsname
\fi
#2}}
\providecommand{\BIBdecl}{\relax}
\BIBdecl

\bibitem{duan2010analysis}
G.-R. Duan, \emph{Analysis and design of descriptor linear systems}.\hskip 1em
  plus 0.5em minus 0.4em\relax Springer Science \& Business Media, 2010,
  vol.~23.

\bibitem{datta2014computation}
S.~Datta and V.~Mehrmann, ``Computation of state reachable points of descriptor
  systems,'' in \emph{Decision and Control (CDC), 2014 IEEE 53rd Annual
  Conference on}.\hskip 1em plus 0.5em minus 0.4em\relax IEEE, 2014, pp.
  6389--6394.

\bibitem{wicks1991computing}
M.~Wicks and R.~A. DeCarlo, ``Computing the distance to an uncontrollable
  system,'' \emph{IEEE Transactions on Automatic Control}, vol.~36, no.~1, pp.
  39--49, 1991.

\bibitem{hu2004real}
G.~Hu and E.~J. Davison, ``Real controllability/stabilizability radius of {LTI}
  systems,'' \emph{IEEE Transactions on Automatic Control}, vol.~49, no.~2, pp.
  254--257, 2004.

\bibitem{gu2006fast}
M.~Gu, E.~Mengi, M.~L. Overton, J.~Xia, and J.~Zhu, ``Fast methods for
  estimating the distance to uncontrollability,'' \emph{SIAM journal on matrix
  analysis and applications}, vol.~28, no.~2, pp. 477--502, 2006.

\bibitem{khare}
S.~R. Khare, H.~K. Pillai, and M.~N. Belur, ``Computing the radius of
  controllability for state space systems,'' \emph{Systems \& Control Letters},
  vol.~61, no.~2, pp. 327--333, 2012.

\bibitem{son}
N.~K. Son and D.~D. Thuan, ``The structured controllability radii of higher
  order systems,'' \emph{Linear Algebra and its Applications}, vol. 438, no.~6,
  pp. 2701--2716, 2013.

\bibitem{son2012structured}
------, ``The structured distance to non-surjectivity and its application to
  calculating the controllability radius of descriptor systems,'' \emph{Journal
  of Mathematical Analysis and Applications}, vol. 388, no.~1, pp. 272--281,
  2012.

\bibitem{simon}
S.~Lam and E.~J. Davison, ``Computation of the real controllability radius and
  minimum-norm perturbations of higher-order, descriptor, and time-delay {LTI}
  systems,'' \emph{IEEE Transactions on Automatic Control}, vol.~59, no.~8, pp.
  2189--2195, 2014.

\bibitem{esing1}
R.~Eising, ``Between controllable and uncontrollable,'' \emph{Systems \&
  control letters}, vol.~4, no.~5, pp. 263--264, 1984.

\bibitem{boley1986measuring}
D.~Boley and W.-S. Lu, ``Measuring how far a controllable system is from an
  uncontrollable one,'' \emph{IEEE Transactions on Automatic Control}, vol.~31,
  no.~3, pp. 249--251, 1986.

\bibitem{elsner1991algorithm}
L.~Elsner and C.~He, ``An algorithm for computing the distance to
  uncontrollability,'' \emph{Systems \& control letters}, vol.~17, no.~6, pp.
  453--464, 1991.

\bibitem{dai}
L.~Dai, \emph{Singular control systems (Lecture notes in control and
  information sciences)}.\hskip 1em plus 0.5em minus 0.4em\relax Springer,
  1989.

\bibitem{lecture}
L.~Scholz, ``Control theory of descriptor systems lecture notes,'' 2015.

\bibitem{Reinschke}
K.~J. Reinschke and G.~Wiedemann, ``Digraph characterization of structural
  controllability for linear descriptor systems,'' \emph{Linear algebra and its
  applications}, vol. 266, pp. 199--217, 1997.

\bibitem{kailath1980linear}
T.~Kailath, \emph{Linear systems}.\hskip 1em plus 0.5em minus 0.4em\relax
  Prentice-Hall Englewood Cliffs, NJ, 1980, vol. 156.

\bibitem{s7}
H.~Park, L.~Zhang, and J.~B. Rosen, ``Low rank approximation of a {H}ankel
  matrix by structured total least norm,'' \emph{BIT Numerical Mathematics},
  vol.~39, no.~4, pp. 757--779, 1999.

\bibitem{s4}
S.~Khare, H.~Pillai, and M.~Belur, ``Numerical algorithm for structured low
  rank approximation problem,'' in \emph{Proceedings of the 19th International
  Symposium on Mathematical Theory of Networks and Systems--MTNS}, vol.~5,
  no.~9, 2010.

\bibitem{s10}
I.~Markovsky, ``Structured low-rank approximation and its applications,''
  \emph{Automatica}, vol.~44, no.~4, pp. 891--909, 2008.

\bibitem{rosen}
J.~B. Rosen, H.~Park, and J.~Glick, ``Total least norm formulation and solution
  for structured problems,'' \emph{SIAM Journal on Matrix Analysis and
  Applications}, vol.~17, no.~1, pp. 110--126, 1996.

\bibitem{gillard2013optimization}
J.~Gillard and A.~Zhigljavsky, ``Optimization challenges in the structured low
  rank approximation problem,'' \emph{Journal of Global Optimization}, vol.~57,
  no.~3, pp. 733--751, 2013.

\bibitem{mengi2008estimation}
E.~Mengi, ``On the estimation of the distance to uncontrollability for higher
  order systems,'' \emph{SIAM Journal on Matrix Analysis and Applications},
  vol.~30, no.~1, pp. 154--172, 2008.

\end{thebibliography}
\end{document}